\documentclass[12pt]{article}
\usepackage{epsfig, epsf, graphicx, subfigure}
\usepackage{pstricks, pst-node, psfrag}
\usepackage{amssymb,amsmath,bm}
\usepackage{verbatim,enumerate}
\usepackage{rotating, lscape}
\usepackage{setspace}

\usepackage{amsfonts,arydshln,breqn}
\usepackage{natbib}
\usepackage{amsfonts,multirow} 

\def\boxit#1{\vbox{\hrule\hbox{\vrule\kern6pt
          \vbox{\kern6pt#1\kern6pt}\kern6pt\vrule}\hrule}}

\def\bse{\begin{eqnarray*}}
\def\ese{\end{eqnarray*}}
\def\be{\begin{eqnarray}}
\def\ee{\end{eqnarray}}
\def\bq{\begin{equation}}
\def\eq{\end{equation}}
\def\bse{\begin{eqnarray*}}
\def\ese{\end{eqnarray*}}

\newcommand{\T}{^\top}
\newcommand{\inv}{^{-1}}
\newcommand{\Real}{\mathbb{R}}
\let\eps=\varepsilon
\let\phi=\varphi

\renewcommand{\red}{\strut}

\begin{document}

\thispagestyle{empty} \baselineskip=28pt \vskip 5mm
\begin{center} {\Huge{\bf On nomenclature for, and the relative merits of, 
   two formulations of skew distributions}}
\end{center}

\baselineskip=12pt \vskip 10mm

\begin{center}\large
Adelchi Azzalini\footnote[1]{
\baselineskip=10pt
Department of Statistical Sciences, University of Padua, 35121 Padova, Italy.
},
Ryan P. Browne\footnote[2]{
\baselineskip=10pt
\red
Department of Statistics and Actuarial Science, University of Waterloo, 
Waterloo, Ontario, Canada, N2L 3G1.
},
Marc G. Genton\footnote[3]{
\baselineskip=10pt 
CEMSE Division,
King Abdullah University of Science and Technology,
Thuwal 23955-6900, Saudi Arabia.
}
and Paul~D.~McNicholas\footnote[4]{
\baselineskip=10pt
Department of Mathematics and Statistics, McMaster University, Hamilton, 
Ontario, Canada, L8S 4L8. 
}

\end{center}

\baselineskip=17pt \vskip 10mm \centerline{\today} \vskip 15mm

%%%%%%%%%%%%%%%%%%%%%%%%%%%%%%%%%%%%%%%%%%%%%%%%%%%%%%%%%%%%%%%%%%%%%%%%
\begin{center}
{\large{\bf Abstract}}
\end{center}
{\red We examine some distributions used extensively within the model-based clustering literature in recent years, paying special attention to} claims that have been made about their relative efficacy. Theoretical arguments are provided as well as real data examples.

\baselineskip=14pt

\par\vfill\noindent
{\bf Some key words:} flexibility; model-based clustering; multivariate distribution; skew-normal distribution; skew-$t$ distribution.
\par\medskip\noindent
{\bf Short title}: Two formulations of skew distributions

\clearpage\pagebreak
%\newpage 
\pagenumbering{arabic}
\baselineskip=26pt

%%%%%%%%%%%%%%%%%%%%%%%%%%%%%%%%%%%%%%%%%%%%%%%%%%%%%%%%%%%%%%%%%%%%%%%%

\section{Introduction} \label{sec:motivation}
{\red In recent years, much work in model-based clustering has replaced the traditional
Gaussian assumption by some more flexible parametric family of distributions. 
In this context,} \cite{leeS:mclachlan:2012sc}, and other work following therefrom, 
utilize two formulations of the multivariate {\red skew-normal} (MSN) distribution as well as analogous formulations of the multivariate skew-$t$ (MST) distribution for clustering, referring to these formulations as ``restricted" and ``unrestricted", respectively. This nomenclature carries obvious implications and, rather than delving into semantics, it will suffice here to quote from \citet[][Section~2.2]{leeS:mclachlan:2012sc}, who contend that ``the unrestricted multivariate skew-normal (uMSN) distribution can be viewed as a simple extension of the rMSN distribution...". Here, rMSN denotes the ``restricted" MSN distribution, and rMST and uMST are used similarly. The purpose of this note is to refute the claim that uMSN distribution is merely a simple extension of the rMSN distribution or, equivalently, the claim that uMST distribution is a simple extension of the rMST distribution. Furthermore, we investigate whether or not one formulation can reasonably be considered superior to the other.

\section{Background}

When one departs from the symmetry of the multivariate normal or other
elliptical distributions, the feature that arises most readily is
skewness. This explains the widespread use of the prefix `skew' which recurs
almost constantly in this context. A recent extensive account is provided by
{\red \citet{azza:capi:2014}}. This activity has generated an enormous number of formulations, sometimes arising with the same motivation and target, or nearly so.  A natural question in these cases is which of the competing alternatives is preferable, either universally or for some given purpose.
To be more specific, start by considering the multivariate skew-normal (SN) distribution proposed by \cite{azza:dval:1996}, examined further by \cite{azza:capi:1999} and by much subsequent work. Note that, although the latter paper adopts a different parameterization of the earlier one, the set of distributions that they
encompass is the same; we shall denote this construction as the classical skew-normal.  Another form of skew-normal distribution has been studied by 
{\red \citet{sahu:dey:branco:2003}}, which we shall refer to as the SDB skew-normal,
{\red by the initials of the author names}. The classical and the SDB set of distributions coincide only for dimension $d=1$; otherwise, the two sets differ and not simply because of different parameterizations. For $d>1$, the question then arises about whether
there is some relevant difference between the two formulations from the viewpoint of suitability for statistical work, both on the side of formal properties and on the side of practical analysis. This question is central to the present note because what we call the classical formulation is what Lee and McLachlan call rMSN, and the SDB formulation is their uMSN.

Analogous formulations arise when the normal family is replaced
by the wider elliptical class in the underlying parent distribution, leading to the so-called skew-elliptical
distributions.  A special case that has received much attention is the
skew-$t$ family \citep{branco:dey:2001,azza:capi:2003}.  Again, the classical skew-$t$ has a
counterpart given by another skew-$t$ considered by
\cite{sahu:dey:branco:2003}, and the same questions as above hold. As before, what we call the classical formulation of the skew-$t$ distribution is what Lee and McLachlan call rMST, and the SDB is their uMST.

Because of their role as the basic constituent for more elaborate 
formulations, we start by discussing the two forms of skew-normal 
distributions. The density and the distribution function of a $N_d(0,\Sigma)$
variable are denoted  $\phi_d(\cdot;\Sigma)$ and  $\Phi_d(\cdot;\Sigma)$,
respectively; the $N(0,1)$ distribution function is denoted 
$\Phi(\cdot)$. The classical skew-normal density function is
\begin{equation}
   f_c(x)= 2\,\phi_d(x-\xi;\Omega)\:\Phi\{\alpha\T\omega\inv(x-\xi)\},
\label{eq:classical-SN-pdf}
\end{equation}
for $x\in\Real^d$, with parameter set $(\xi, \Omega, \alpha)$. Here $\xi$ is a $d$-dimensional 
location parameter, $\Omega$ is a symmetric positive definite $d\times d$ 
scale matrix, $\alpha$ is a 
$d$-dimensional slant parameter, and $\omega$ is a diagonal matrix formed by
the square roots of the diagonal elements of $\Omega$.
Various stochastic representations exist for (\ref{eq:classical-SN-pdf}). 
One is as follows: if
$$   \begin{pmatrix}
   X_0 \\ X_1
   \end{pmatrix}
   \sim N_{d+1}(0, \Omega^*), \qquad 
   \Omega^* = \begin{pmatrix} \bar\Omega   & \delta \\
                        \delta\T      & 1,
   \end{pmatrix}$$
where $\Omega^*$ is a correlation matrix, then
\begin{equation}
   Y_c = \xi + \omega (X_0|X_1>0) 
  \label{eq:classical-SN-sr1}
\end{equation}
has distribution (\ref{eq:classical-SN-pdf}) with
$\Omega=\omega\bar\Omega\omega$ and 
$ \alpha =(1-\delta\T\bar\Omega\inv\delta)^{-1/2}\:\bar\Omega\inv\delta.$
{\red Here and in the following, given a random variable  $X$ 
and  an event $E$, the notation $(X|E)$ denotes a random variable 
which has the distribution of $X$ conditional on the event $E$;
the Kolmogorov representation theorem ensures that such a random variable
exists.}

Another stochastic representation is the following: if $\delta$ is a
$d$-vector with elements in $(-1,1)$, then (\ref{eq:classical-SN-pdf}) is the
density function of
\begin{equation}
  Y_c = \xi + \omega\left\{ 
        [I_d- \mathrm{diag}(\delta)^2]^{1/2}\, V_0 + \delta|V_1|\right\},
  \label{eq:classical-SN-sr2}
\end{equation}
where $V_0$ and $V_1$ are independent normal variates of dimension $d$ and $1$,
respectively, with 0 mean value, unit variances, and $\mathrm{cor}(V_0)$ is
suitably related to $\alpha$ and $\Omega$; {\red  full details are given
on p.\,128--9 of  \citet{azza:capi:2014} among other sources.}
 For the SDB skew-normal, we adopt a very minor change from the symbols of \cite{sahu:dey:branco:2003}, but retain the same parameterization.  
Given real values $\lambda_1, \dots,\lambda_d$, let $\lambda=(\lambda_1, \dots, \lambda_d)\T$ and $\Lambda=\mathrm{diag}(\lambda)$, and write the SDB density as
\begin{dmath}%\begin{split}
  f_s(x) = 2^d\,\phi_d(x-\xi; \Delta+\Lambda^2)
       \times\Phi_d\{\Lambda(\Delta+\Lambda^2)\inv(x-\xi); 
       I_d-\Lambda(\Delta+\Lambda^2)\inv\Lambda\},
\label{eq:SDB-SN-pdf}
%\end{split}
\end{dmath}
where $\Delta$ is a symmetric positive-definite matrix.
This density is associated with the following stochastic representation. 
For independent variables $\eps\sim N_d(\xi, \Delta)$ and 
$Z\sim N_d(0, I_d)$, consider the  transformation
\begin{equation}
     Y_s= \Lambda (Z|Z>0) + \eps,
\label{eq:SDB-SN-sr}
\end{equation}
where $Z>0$ means that the inequality is satisfied component-wise; 
then $Y_s$ has density (\ref{eq:SDB-SN-pdf}).

\section{Comparing the Formulations}

A qualitative comparison of the formal properties of the two distributions lends 
several annotations.  Some of these have already been presented by 
\cite{sahu:dey:branco:2003}, but they are included here for completeness.
\begin{enumerate}
\item The number of individual parameter values is $2d+d(d+1)/2$ in both cases.
\item The two families of distributions coincide only for $d=1$, {\red as noted by 
 \cite{sahu:dey:branco:2003},} and neither one is a subset of the other for $d>1$.
\item {\red As $d$ increases}, computation of $f_s$ becomes {\red progressively more} cumbersome because of 
  {\red the factor} $\Phi_d$. 
\item \label{item:affine} The classical skew-normal family is closed under affine transformations,   while the same fact does not hold for the SDB family.
\item \label{item:factorize} Another remark of \cite{sahu:dey:branco:2003} is that $f_s$ can allow for $d$ independent skew-normal components, when $\Delta$ is diagonal, while $f_c$ can factorize only as a product where at most one factor is skew-normal with non-vanishing slant parameter. 
\item \label{item:latent} Stochastic representations (\ref{eq:classical-SN-sr1}) and (\ref{eq:SDB-SN-sr}) involve $1$ and $d$ latent variables, respectively. The latter one seems to fit less easily in an applied setting, because it requires that for each observed component there is a matching latent component, while the classical construction can more easily be incorporated in the logical frame describing  a real phenomenon subject to selective sampling based on one latent variable.
\item For the classical skew-normal, the expressions of higher order cumulants and Mardia's coefficients of multivariate skewness and kurtosis are given in Appendix A.2 of \citet{azza:capi:1999}. The range of skewness is $[0,  g_1^*)$ where $g_1^*=2(4-\pi)^2/(\pi-2)^3$; the range of excess kurtosis is $[0, g_2^*)$, where $g_2^*=8(\pi-3)/(\pi-2)^2$. For the SDB form, expressions of Mardia's coefficients are given in the Appendix. Numerical maximization of these expressions when $d=2$ leads to ranges with maximal values that appear to coincide numerically with $2\,g_1^*$ and $2\,g_2^*$, respectively.
\item For the classical skew-normal, the distribution of quadratic forms can be obtained from the similar case under normality. No similar result is known to hold for the SDB form.
\end{enumerate}  
Clearly, these remarks do not lead one to consider either one formulation superior to the other.

Each of the two skew-normal families discussed above leads to a matching form
of skew-$t$ family. For the classical case, this can be obtained
by replacing the assumption of joint normality of $(X_0, X_1)$ in
(\ref{eq:classical-SN-sr1}) by one of $(d+1)$-dimensional Student's $t$
distribution \citep{branco:dey:2001,azza:capi:2003}. The SDB skew-$t$ has been obtained by
\cite{sahu:dey:branco:2003} assuming that $(Z,\eps)$ entering
(\ref{eq:SDB-SN-sr}) is a $(2d)$-dimensional Student's~$t$.
In both cases, the resulting density is similar in structure to the
skew-normal case, with the $\phi_d$ {\red factor} replaced by a $d$-dimensional $t$
density on $\nu$ degrees of freedom, but the skewing factor is different: for
the classical version, it is given by the distribution function of a univariate
$t$ on $\nu+d$ degrees of freedom; for the SDB version, the $t$ distribution
function is $d$-dimensional. 

We now return to the relationship between the classical and SDB formulations,
{\red as discussed by \citet{leeS:mclachlan:2012sc}:} ``The unrestricted
multivariate skew-normal (uMSN) distribution can be viewed as a simple
extension of the rMSN distribution in which the univariate latent variable
$U_0$ is replaced by a multivariate analogue, that is, $U_0$.''  Note that
their $U_0$ is our $V_1$ in (\ref{eq:classical-SN-sr2}). In reality, the use
of a multivariate latent error term in place of a single random component does
not add any level of generality because this multivariate latent variable,
essentially $Z\sim N_d(0, I_d)$ in (\ref{eq:SDB-SN-sr}), has a highly
restricted structure. It entails more random ingredients than the single
$X_1\sim N(0,1)$ variable in (\ref{eq:classical-SN-sr1}); however, because of
the highly restricted structure, {\red we are not provided with more parameters 
to maneuver for   increasing flexibility, as already indicated by the
fact that the overall number of parameters is the same in the two formulations}.
Furthermore, the use of ``extension'' is clearly inappropriate
because neither one of the two families is a subset of the other for $d>1$.

{\red A further relevant aspect appears} in the discussion of the classical skew-$t$
distribution near the end of Section~3 of \cite{leeS:mclachlan:2012sc},  {\red where}
the authors state that ``the form of skewness is limited in these
characterizations. In Sect.\,5 we study an extension of their approach to the
more general form of skew $t$-density as proposed by Sahu et al. (2003).''
This claim of limited form of skewness is supported only by the
above-indicated misinterpretation of the role of the perturbation factor, and
not by any concrete elements. 
Indeed, if one looks at quantitative elements, the opposite message emerges, {\red 
as explained next. First of all, note that the wider ranges of the Mardia's
measures of skewness and kurtosis for the SDB skew-normal distribution, as
mentioned earlier in this section, are of little relevance because the range of these
measures is very limited anyway. To achieve a substantial level of skewness
and kurtosis one has to adopt some form of skew-$t$ distribution with small degrees of
freedom. Now,}
the range of skewness for the classical skew-$t$ distribution is unlimited both marginally, when measured by the usual coefficient $\gamma_1$, and globally, when measured by Mardia's coefficient $\gamma_{1,d}$. To see this fact in the univariate case, which coincides with the behaviour of univariate components in a multivariate skew-$t$ distribution, consider the expression of $\gamma_1$ on p.\,382 of \citet{azza:capi:2003} and let $\nu\to3$; for the Mardia's coefficient, see (6.31) on p.\,178 of \citet{azza:capi:2014}.

\section{Model-Based Clustering Illustrations} 

Because model-based clustering represents the context in which the
nomenclature under consideration has been popularized {\red },
illustrations will be focused in that direction. 
In brief, model-based clustering is the use of (finite) mixture
models for clustering. A finite mixture of rMST (FM-rMST) distributions is
simply a convex linear combination of rMST distributions, and FM-uMST has an
analogous meaning. \cite{leeS:mclachlan:2012sc} provide numerous clustering
illustrations where the ``superiority" and ``extra flexibility" of the FM-uMST
are illustrated. While it is true that such terminology is used 
{\red } in relevance to particular data sets or examples,
it is also true that there are many such examples and all have more or less
the same message: the FM-uMST is better, in some sense, than the FM-rMST. This
pattern is also present in other work by the same authors, {\red for instance in}
\citet{leeSX:mclachlan:2013adac}. The goal of the analyses herein is
to present an extensive comparison using algorithms written by \citet{skew}
and by \citet{uskew}.

To avoid the perception of bias that can arise from selection of a subset of variables from a given data set, we examine \emph{all possible} pairs and triplets of variables in two real data sets that are commonly used in model-based clustering illustrations. Although these illustrations are conducted as genuine cluster analyses, i.e., without knowledge of labels, the labels are known; accordingly, we can assess the classification performance of the fitted models. The adjusted Rand index \citep[ARI;][]{hubert85} is used for this purpose. It takes a value of one when there is perfect agreement between two classes, and its expected value is zero under random classification. 
We consider the crabs data \citep{campbell74} {\red available in} the {\tt
  MASS} package for {\sf R} {\red \citep{Rpkg-MASS,R14}}. These data comprise five biological measurements of 200 crabs of genus \textit{Leptograpus}, 
  i.e., 50 male and 50 female crabs for each of two species. 
We also consider the Australian Institute of Sport (AIS) data, which comprise 11 biomedical and anthropometric measurements and two categorical variables, i.e., gender and sport, for each of 202 Australian athletes. 
For each data set, we proceed  by considering all possible pairs and triplets of the continuous measurements to build  clusters of the data points --- a total of 480 cluster analyses.  
{\red As is standard practice \citep[e.g.][]{peel:mclachlan:2000},}
we take gender as the reference label when computing the ARI.
The \textsf{R} packages \texttt{EMMIXskew} \citep{skew} and \texttt{EMMIXuskew} \citep{uskew} are used to implement the ``restricted" and ``unrestricted" formulations, respectively, and the {\sf R} code we use to produce the results in this section is available in Supplementary Material.
The results (Figure~\ref{fig:comp}) very clearly indicate that neither formulation is markedly  superior and, if these results were to be taken in favour of either formulation, it would be the classical formulation.
%\begin{figure}[!bt]
%\centering
%\vspace{-0.25in}
%  \includegraphics[width=0.92\hsize]{crabs56} \\
%\vspace{-0.12in}
%  \includegraphics[width=0.92\hsize]{ais56}
%%\vspace{-0.2in}
%  \caption{ARI values from model-based clustering analyses of the crabs (top) and AIS (bottom) data using the \texttt{EMMIXskew} and \texttt{EMMIXuskew} packages, where pairs are represented by circles and triplets by triangles.}% eps=1e-3 and itmax=1000
%  \label{fig:comp}
%\end{figure}
\begin{figure*}[!bt]
\centering
%\hspace{-0.1in}
  \includegraphics[width=0.495\hsize]{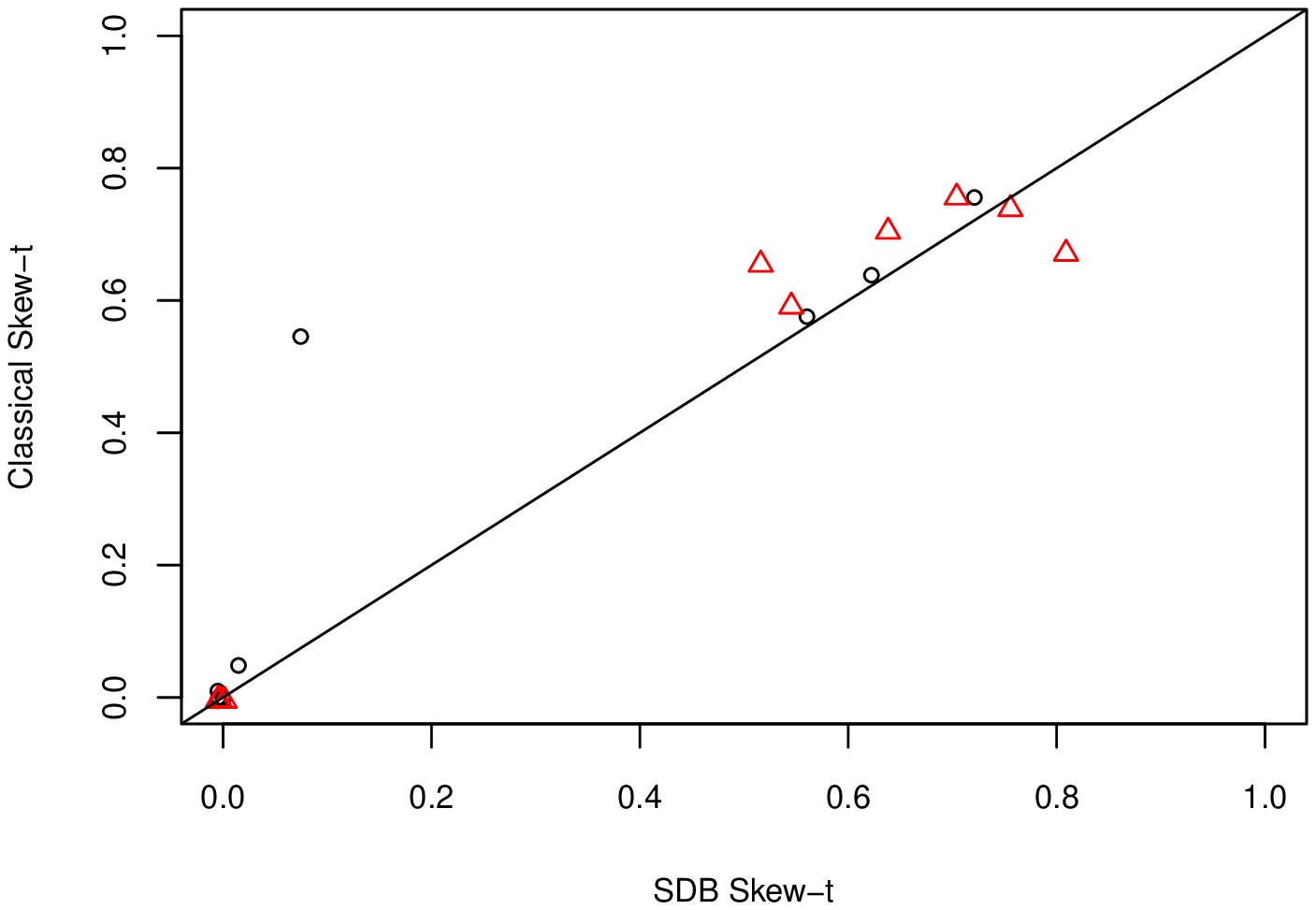} 
  \includegraphics[width=0.495\hsize]{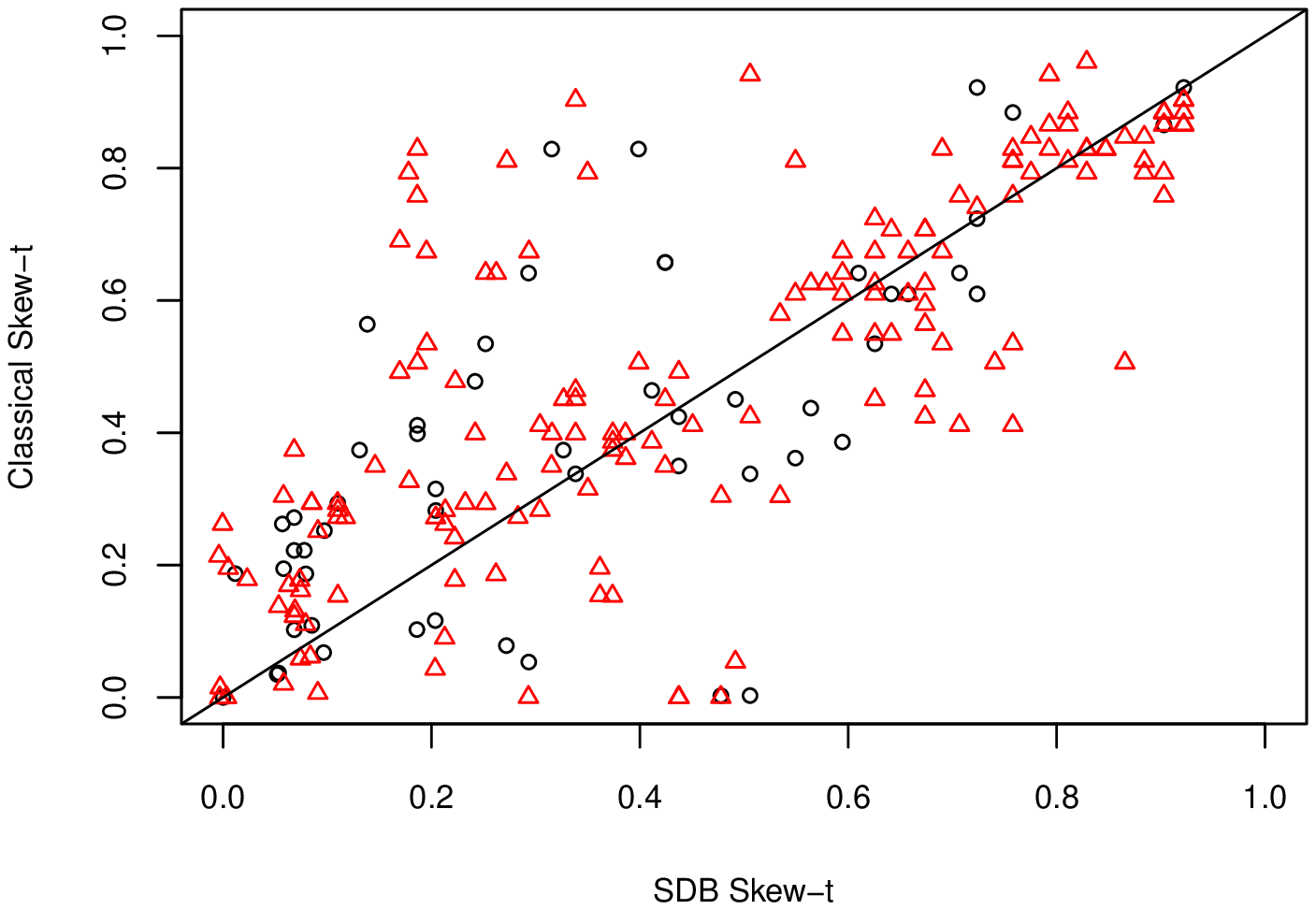}
\vspace{-0.2in}
 \caption{ARI values from model-based clustering analyses of the crabs (left) and AIS (right) data using the \texttt{EMMIXskew} and \texttt{EMMIXuskew} packages, where pairs are represented by circles and triplets by triangles.}% eps=1e-3 and itmax=1000
  \label{fig:comp}
\end{figure*}

In addition to the results in Figure~\ref{fig:comp}, which correspond to {\tt
  SEED=5}, the Supplementary Material can be used to easily produce results
from other starting values by modifying {\tt SEED}. As the reader can verify,
running the code for {\tt SEED=1,...,10} leads to the same message, i.e.,
neither formulation is superior, but now based on 4,800 cluster analyses. It
is noteworthy {\red  that the only scenario where the SDB formulation produced better clustering results was constructed by starting the numerical optimization
  using the true classification labels as the initial values in the numerical
  search}.  Of course, assessing a clustering method based on how it does
when started at the true classifications is fundamentally flawed because {\red
  the true classifications are not available} in real cluster analyses,
{\red assuming that a true classification even exists}.
 Furthermore,
such an assessment will lead to an excellent assessment for any method that
just does not move from the starts --- even a method that simply returns the
starting classifications.

A final matter for consideration is the relative computation time for the two
formulations.  
{\red With few exceptions, examples within the literature where the SDB 
formulation is used for meaningful analysis only consider data with  $d\leq 4$;}
see e.\,g.\ \citet[Section~6]{leeS:mclachlan:2012sc}. 
It is instructive to consider
the ratio of time taken by the SDB formulation to the time taken by the
classical formulation for the pairs and triples of the AIS and crabs data. The
results of this comparison (Table~\ref{tab:ratio}) confirm that the SDB
formulation is very much slower than the classical formulation, e.g., taking
an average of 7,315 times longer to converge for the three-dimensional crabs
data. The {\sf R} code used to produce the results in Table~\ref{tab:ratio} is
available in Supplementary Material.
% latex table generated in R 3.1.2 by xtable 1.7-3 package
% Thu Apr  9 09:20:37 2015
\begin{table}[t]
\centering
\caption{Means and standard deviations for the ratios of {\tt user} times,
  from the {\tt system.time()} function in {\sf R}, for the SDB formulation to
  the classical formulation when applied to two- and three-dimensional subsets
  of the AIS and crabs data.\label{tab:ratio}}
\par\vspace{2ex}
\begin{tabular*}{0.49\textwidth}{@{\extracolsep{\fill}}lrrr}
 \hline
&Dimension& Mean & Std.\ Deviation \\ 
 \hline
\multirow{2}{*}{AIS} & 2 & 125.4 & 160.4 \\ 
 & 3 & 2216.0 & 2131.9 \\ 
\hline
 \multirow{2}{*}{Crabs} & 2 & 307.5 & 93.5 \\ 
 & 3 & 7315.1 & 3619.2 \\ 
  \hline
\end{tabular*}
\end{table}

\section{Conclusion}\label{sec:conc}

{\red
We have discussed the relative merits of two closely related formulations, each leading 
to a skewed extension of the multivariate normal and  $t$ distribution.
 For one, we have clarified} why the SDB (or ``unrestricted'') formulation is
not a ``simple extension" of the classical (or ``restricted") formulation. We
also provide extensive evidence as to why neither formulation is, in general,
preferable to the other. Extensive {\red numerical work (4,800 cases in all) was} carried out to underline this point in specific reference to clustering applications. We trust it is now clear that the nomenclature ``restricted" and ``unrestricted" should be avoided in reference to these formulations. 

\section*{Acknowledgements}
We are grateful to M\'arcia Branco for fruitful discussions of various aspects
of the SDB formulation and for making available to us related material.

\section*{Appendix}

{\bf Cumulants and Mardia's coefficients for SDB skew-normal}. 
For the SDB skew-normal,  differentiation of $K_s(t)=\log M_s(t)$ produces
\begin{eqnarray*}
   \nabla K_s(t) &=& \xi + (\Delta+\Lambda^2) t + 
        \left[\zeta_1(\lambda_j t_j)\lambda_j\right]_{j=1}^d   \label{e:K'}\\
  \nabla\nabla\T K_s(t) &=& (\Delta+\Lambda^2) + \mathrm{diag}(\zeta_2(\lambda_1 t_1) \lambda_1^2,\dots, 
         \zeta_2(\lambda_d t_d) \lambda_d^2),  \label{e:K"}
\end{eqnarray*}
where $\zeta_r(x)$ is the $r$th derivatives of
$\zeta_0(x)=\log\{2\,\Phi(x)\}$.  Evaluation at $t=0$ gives
%(\ref{eq:SDB-SN-mean,var})
\begin{equation}
  E(Y_s) =\xi + \sqrt{2/\pi}\: \lambda, \ \ 
  \mathrm{var}(Y_s)=\Delta + (1-2/\pi) \Lambda^2.
\label{eq:SDB-SN-mean,var}
\end{equation}

Further differentiation and evaluation at $0$ gives the 3rd order cumulant
\begin{equation}\label{e:kappa3}
   \kappa_{rst} = \left\{ \begin{array}{ll} \zeta_3(0)\, \lambda_r^3 & \mbox{if } r=s=t,\\
                         0           & \mbox{otherwise}, \end{array} \right.
\end{equation}
where $\zeta_3(0)=b\,(4/\pi-1) = (2/\pi)^{3/2}\:(4-\pi)/2$.

We can now compute the Mardia's coefficient $\gamma_{1,d}$ of multivariate
skewness, recalling that the
3rd order cumulant coincides with the 3rd order central moment.  Denote by
$\Sigma=(\sigma_{rs})$ the variance matrix in (\ref{eq:SDB-SN-mean,var}) and
let $\Sigma\inv= (\sigma^{rs})$, $\mu_j=b\lambda_j$.

From (2.19) of \citet{mardia:1970}, write
\begin{eqnarray}
  \gamma_{1,d} &=& \sum_{rst} \sum_{r's't'} \kappa_{rst} \kappa_{r's't'}
                    \sigma^{rr'}  \sigma^{ss'}  \sigma^{tt'} 
       = \zeta_3(0)^2 \sum_{u,v} \lambda_u^3\,\lambda_v^3\,(\sigma^{uv})^3
            \nonumber\\
       &=& \left(\frac{4-\pi}{2}\right)^2 
                \sum_{u,v} \mu_u^3\,\mu_v^3\,(\sigma^{uv})^3  
       = \left(\frac{4-\pi}{2}\right)^2 (\mu^{(3)})\T\Sigma^{(-3)}\mu^{(3)},
           \label{e:gamma1M-altern}
\end{eqnarray}
where $\mu^{(3)}$ is the vector with elements $\mu_j^3$ and
$\Sigma^{(-3)}=\left((\sigma^{uv})^3\right)$.

For (\ref{e:gamma1M-altern}) we do not have an expression of the maximal value.
Numerical exploration indicates that the maximal value is  $1.98113$, that is, 
the double  value of the classical SN {\red up to the quoted number of digits}.
This maximal value of
$\gamma_{1,d}$ is obtained, irrespectively of $\Delta$, in these four cases:
\begin{equation}
   \lambda  = h (\pm 1 \,\, \pm1)\T, \quad \hbox{when~} h\to\infty \,.
  \label{e:extreme-lambda}
\end{equation}
%	A graphical illustration of the corresponding type of density is
% provided by the graph below, which refers to $\Delta=I_2$,
% $\lambda=10^3(1,1)\T$.  The other three cases of $\lambda$ produce similar
% plots rotated by $\pm\pi/2$ or $\pi$, placed in the other three quadrants of
% $\Real^2$.\par
%	\centerline{\includegraphics[width=0.5\hsize]{Figures/pdf-SDB-1000}}

Derivation of the 4th order cumulants is similar to (\ref{e:kappa3}), leading
to
\begin{equation} \label{e:kappa4}
   \kappa_{rstu} = \left\{ \begin{array}{ll} \zeta_4(0)\, \lambda_r^4 & \mbox{if } r=s=t=u,\\
                         0           & \mbox{otherwise}, \end{array} \right.
\end{equation}
where $\zeta_4(0)=  2\:(\pi-3)\:(2/\pi)^2 \approx 0.114771$. 

From here the Mardia's coefficient of (excess) kurtosis is
\begin{eqnarray}
  \gamma_{2,d} &=& \sum_{rstu} \kappa_{rstu} \sigma^{rs}  \sigma^{tu} 
             \nonumber  \\
       &=& \zeta_4(0) \sum_{u} \lambda_u^4\,(\sigma^{uu})^2  \nonumber\\
       &=& 2\,(\pi-3) \sum_{u} \mu_u^4\,\,(\sigma^{uu})^2  \nonumber\\
       &=& 2\,(\pi-3)  (\mu^{(2)})\T(I_d \odot \Sigma\inv)^2\mu^{(2)},
           \label{e:gamma2M-altern}
\end{eqnarray}
where $\mu^{(2)}=(\mu_1^2,\dots, \mu_d^2)\T$ and $\odot$ is the Hadamard
{\red or component-wise} product.  
A numerical search indicates that the maximal value of
$\gamma_{2,d}$ is again achieved, irrespectively of $\Delta$, with $\lambda$
as in (\ref{e:extreme-lambda}).  The maximal observed value of the coefficient
is $1.7383546$, again twice the corresponding value in the classical 
skew-normal case.


{\small \begin{thebibliography}{}

\bibitem[\protect\citeauthoryear{Azzalini and Capitanio}{Azzalini and
  Capitanio}{1999}]{azza:capi:1999}
Azzalini, A. and A.~Capitanio (1999).
\newblock Statistical applications of the multivariate skew normal
  distribution.
\newblock {\em J. Roy. Stat. Soc., series B\/}~{\em 61\/}(3), 579--602.

\bibitem[\protect\citeauthoryear{Azzalini and Capitanio}{Azzalini and
  Capitanio}{2003}]{azza:capi:2003}
Azzalini, A. and A.~Capitanio (2003). 
\newblock Distributions generated by perturbation of symmetry with emphasis on
  a multivariate skew $t$ distribution.
\newblock {\em J. Roy. Stat. Soc., series B\/}~{\em 65\/}(2), 367--389.

\bibitem[\protect\citeauthoryear{Azzalini}{Azzalini and Capitanio}{2014}]{azza:capi:2014}
Azzalini, A. with the collaboration of A.~Capitanio (2014).
\newblock {\em The Skew-Normal and Related Families}.
\newblock IMS monographs. Cambridge University Press.

\bibitem[\protect\citeauthoryear{Azzalini and {Dalla Valle}}{Azzalini and
  {Dalla Valle}}{1996}]{azza:dval:1996}
Azzalini, A. and A.~{Dalla Valle} (1996).
\newblock The multivariate skew-normal distribution.
\newblock {\em Biometrika\/}~{\em 83}, 715--726.

\bibitem[\protect\citeauthoryear{Branco and Dey}{Branco and
  Dey}{2001}]{branco:dey:2001}
Branco, M.~D. and D.~K. Dey (2001).
\newblock A general class of multivariate skew-elliptical distributions.
\newblock {\em J. Multivariate Analysis\/}~{\em 79}, 99--113.

\bibitem[\protect\citeauthoryear{Campbell and Mahon}{Campbell and
  Mahon}{1974}]{campbell74}
Campbell, N.~A. and R.~J. Mahon (1974).
\newblock A multivariate study of variation in two species of rock crab of
  genus leptograpsus.
\newblock {\em Australian Journal of Zoology\/}~{\em 22}, 417--425.

\bibitem[\protect\citeauthoryear{Hubert and Arabie}{Hubert and
  Arabie}{1985}]{hubert85}
Hubert, L. and P.~Arabie (1985).
\newblock Comparing partitions.
\newblock {\em Journal of Classification\/}~{\em 2}, 193--218.

\bibitem[\protect\citeauthoryear{Lee and {McLachlan}}{Lee and
  {McLachlan}}{2014}]{leeS:mclachlan:2012sc}
Lee, S. and G.~J. {McLachlan} (2014).
\newblock Finite mixtures of multivariate skew $t$-distributions: some recent
  and new results.
\newblock {\em Statistics and Computing\/}~{\em 24}, 181--202.
\newblock Published online 20 October 2012.

\bibitem[\protect\citeauthoryear{Lee and McLachlan}{Lee and
  McLachlan}{2013a}]{leeSX:mclachlan:2013adac}
Lee, S.~X. and G.~J. McLachlan (2013a).
\newblock On mixtures of skew normal and skew $t$-distributions.
\newblock {\em Advances in Data Analysis and Classification\/}~{\em 7},
  241--266.
  
\bibitem[\protect\citeauthoryear{Lee and McLachlan}{Lee and
  McLachlan}{2013b}]{uskew}
Lee, S.~X. and G.~J. McLachlan (2013b).
\newblock {EMMIXuskew}: An {R} package for fitting mixtures of multivariate
  skew $t$ distributions via the {EM} algorithm.
\newblock {\em Journal of Statistical Software\/}~{\em 55\/}(12), 1--22.


%\bibitem[\protect\citeauthoryear{Lee and {McLachlan}}{Lee and
%  {McLachlan}}{2015}]{leeSX:mclachlan:2014[8182]}
%Lee, S.~X. and G.~J. {McLachlan} (2015).
%\newblock Finite mixtures of canonical fundamental skew $t$-distributions.
%\newblock {\em Statistics and Computing\/}. To appear.
%\newblock Published online 28 February 2015.

\bibitem[\protect\citeauthoryear{Mardia}{Mardia}{1970}]{mardia:1970}
Mardia, K. (1970).
\newblock Measures of multivariate skewness and kurtosis with applications.
\newblock {\em Biometrika\/}~{\em 57}, 519--530.


\bibitem[\protect\citeauthoryear{Peel and {McLachlan}}{Peel and
  {McLachlan}}{2000}]{peel:mclachlan:2000}
Peel, D. and G.~J. {McLachlan} (2000).
\newblock  Robust mixture modelling using the $t$ distribution.
\newblock {\em Statistics and Computing\/}~{\em 10}, 339--348.

\bibitem[\protect\citeauthoryear{R Core Team}{R Core Team}{2014}]{R14}
R Core Team (2014).
\newblock R: A Language and Environment for Statistical Computing,
\newblock R Foundation for Statistical Computing,Vienna, Austria. % Version 3.1.2. 
%(2014-10-31) -- "Pumpkin Helmet"

\bibitem[\protect\citeauthoryear{Sahu, Dey, and Branco}{Sahu
  et~al.}{2003}]{sahu:dey:branco:2003}
Sahu, K., D.~K. Dey, and M.~D. Branco (2003).
\newblock A new class of multivariate skew distributions with applications to
  {B}ayesian regression models.
\newblock {\em Canad. J. Statist.\/}~{\em 31\/}(2), 129--150.
\newblock Corrigendum: vol.\ 37 (2009), 301--302.

\bibitem[Venables \& Ripley, 2002]{Rpkg-MASS}
{\red  Venables, W.~N. \& Ripley, B.~D. (2002).
  \newblock {\em Modern Applied Statistics with S}.
  \newblock New York: Springer, fourth edition.}
  % \newblock ISBN 0-387-95457-0.

\bibitem[\protect\citeauthoryear{Wang, Ng, and McLachlan.}{Wang
  et~al.}{2013}]{skew}
Wang, K., A.~Ng, and G.~McLachlan. (2013).
\newblock {\em EMMIXskew: The EM Algorithm and Skew Mixture Distribution}.
\newblock R package version 1.0.1.

\end{thebibliography}}
\end{document}